\documentclass[onefignum,onetabnum]{siamart171218}

\usepackage{lipsum}
\usepackage{amsfonts}
\usepackage{graphicx}
\usepackage{epstopdf}
\usepackage{algorithmic}
\ifpdf
  \DeclareGraphicsExtensions{.eps,.pdf,.png,.jpg}
\else
  \DeclareGraphicsExtensions{.eps}
\fi

\newsiamremark{remark}{Remark}
\newsiamremark{hypothesis}{Hypothesis}
\newsiamremark{conjecture}{Conjecture}
\crefname{hypothesis}{Hypothesis}{Hypotheses}
\newsiamthm{claim}{Claim}

\headers{Numerical Conformal Mappings}{H. Hakula, A. Rasila}

\title{Laplace--Beltrami Equations and Numerical Conformal Mappings on Surfaces}

\author{Harri Hakula\thanks{Department of Mathematics and System Analysis, Aalto University, P.O. Box 11100,
FI--00076 Aalto, Finland 
  (\email{Harri.Hakula@aalto.fi}).}
\and Antti Rasila\thanks{Department of Mathematics with Computer Science, Guangdong Technion -- Israel Institute of Technology, Shantou, Guangdong 515063, P.R. of China and Department of Mathematics with Computer Science, Technion -- Israel Institute of Technology
Haifa 32000, Israel. Corresponding author.
  (\email{antti.rasila@iki.fi}, \email{antti.rasila@gtiit.edu.cn}).}
}

\usepackage{amsopn}

\makeatletter
\newcommand*{\addFileDependency}[1]{
  \typeout{(#1)}
  \@addtofilelist{#1}
  \IfFileExists{#1}{}{\typeout{No file #1.}}
}
\makeatother

\usepackage{subfig}
\usepackage{interval}
\intervalconfig{soft open fences}

\usepackage{booktabs}
\usepackage{enumitem}

\newcommand{\cT}{\mathcal{T}}
\newcommand{\eE}{\varepsilon}
\newcommand{\osc}{\mathrm{osc}}
\newcommand{\vr}{\mathbf{r}} 
\newcommand{\Sf}{\mathcal{S}} 

\newcommand{\C}{\mathbb{C}}

\graphicspath{{PDF/}{Graphics/}}

\ifpdf
\hypersetup{
  pdftitle={Laplace--Beltrami Equations and Numerical Conformal Mappings on Surfaces},
  pdfauthor={H. Hakula, A. Rasila}
}
\fi

\begin{document}

\maketitle

\begin{abstract}
  The conjugate function method is an algorithm for numerical computation of conformal mappings for simply and multiply connected domains.
    In this paper, the conjugate function method is extended
    to cover conformal mappings between Riemannian surfaces. 
    The main challenge addressed here is the connection
    between Laplace--Beltrami equations on surfaces and the computation of the conformal modulus of a quadrilateral.
    We consider mappings of simply, doubly, and multiply connected domains. The numerical computation is based on an $hp$-adaptive finite element method. The key advantage of our approach is that it allows highly accurate computations of mappings on surfaces, including domains of complex boundary geometry involving strong singularities and cusps.
    The efficacy of the proposed method is illustrated via an extensive set of numerical experiments including error estimates.
\end{abstract}



\section{Introduction}
\label{sec:intro}

The problem of finding numerical conformal mappings between plane domains has been widely studied in the literature. Popular methods in this setting include the Schwarz-Christoffel Mapping \cite{Driscoll,driscoll_trefethen_2002}, the Zipper algorithm \cite{Marshall}, and boundary integral methods \cite{Nasser,Wegmann-Nasser}. Our earlier work in this topic includes \cite{hqr}, where the conjugate function method was introduced and its generalization to multiply connected domains in \cite{hqr2}. Both of these papers make use of the $hp$-FEM method introduced in \cite{hrv} for numerical computation of underlying conformal moduli and potential functions required by the method, but other approaches can be used as well. For example in \cite{hrs}, conformal moduli and potential functions are approximated by using a stochastic Walk-on-Spheres method. For comprehensive surveys on methods available for numerical conformal mappings in the plane, see e.g. \cite[pg. 8--11]{Kythe} and \cite{Wegmann}.

While, as in the planar case, the existence of conformal mappings between simply connected Riemannian surfaces is well-known, the possibility of applying aforementioned algorithms in this setting is not clear, making this topic substantially more difficult than obtaining mappings between plane domains. Numerical methods for constructing conformal mappings between surfaces include circle packings \cite{COLLINS2003233,Sullivan,Thurston} and the differential geometric approach of Gu and Yau \cite{Gu-Yau}.


In this paper, we generalize our earlier work on numerical conformal mappings between plane domains to mappings of a canonical plane domain onto a surface of the same topological type. Our method is based on conformal modulus, which is applicable in very general settings, and on finding harmonic solutions to Dirichlet-Neumann mixed boundary value problems on surfaces through state-of-the-art numerical methods. As a result, we obtain an algorithm that is not only easy to understand but also very general and accurate even for surfaces of complicated geometric types.

\subsection{PDEs on Surfaces}
Partial differential equations defined on general surfaces arise in many different contexts.
Several methods for defining suitable discretizations of surfaces and corresponding 
finite element spaces have been proposed. In this work, it is assumed that 
some global parameterization is available for each instance.
It should be mentioned that 
Demlow~\cite{demlow2007adaptive} mentions that this approach is not sufficiently flexible to cover many
interesting applications such as those with evolving or moving surfaces, indeed,
such problems are not considered here.
There are many excellent references covering approaches where the surface
is approximated by a polyhedral surface having triangular faces, and we only list a few, see  for instance
\cite{demlow2012adaptive,bonito2020finite}.
Demlow also considers higher order discretizations within the classical $h$-version of the finite element method \cite{demlow2009higher} and 
Cantwell et al. \cite{CANTWELL2014813} within the $hp$-version.
Bonito and Demlow have contributed an implementation to a popular open source library \texttt{deal.II} \cite{dealii2019design}.
Also, some recent work covers more modern discretization techniques such as virtual finite elements \cite{refId0}.

\subsection{Conformal Mappings on Surfaces}

The following counterpart of the Riemann mapping theorem for surfaces is known as the uniformization theorem (see e.g. \cite[Theorem 10-3]{Ahlfors}):


\begin{theorem}
Every simply connected Riemann surface is conformally equivalent to a disk, to the complex plane, or to the Riemann sphere.
\end{theorem}

This result was independently proved  by Paul Koebe and Henri Poincar\'e in 1907.

In this paper, we also investigate conformal mappings of doubly and multiply connected surfaces, where the mapping is constructed onto canonical domains as in \cite{hqr,hqr2}. In these cases, the existence of the conformal mapping is guaranteed by a more general version of the above result, see e.g. \cite[Theorem 3.1]{Jost}.

\subsection{Illustrative Example: Cartography}
\begin{figure}
  \centering
  \subfloat[]{
    \includegraphics[height=2in]{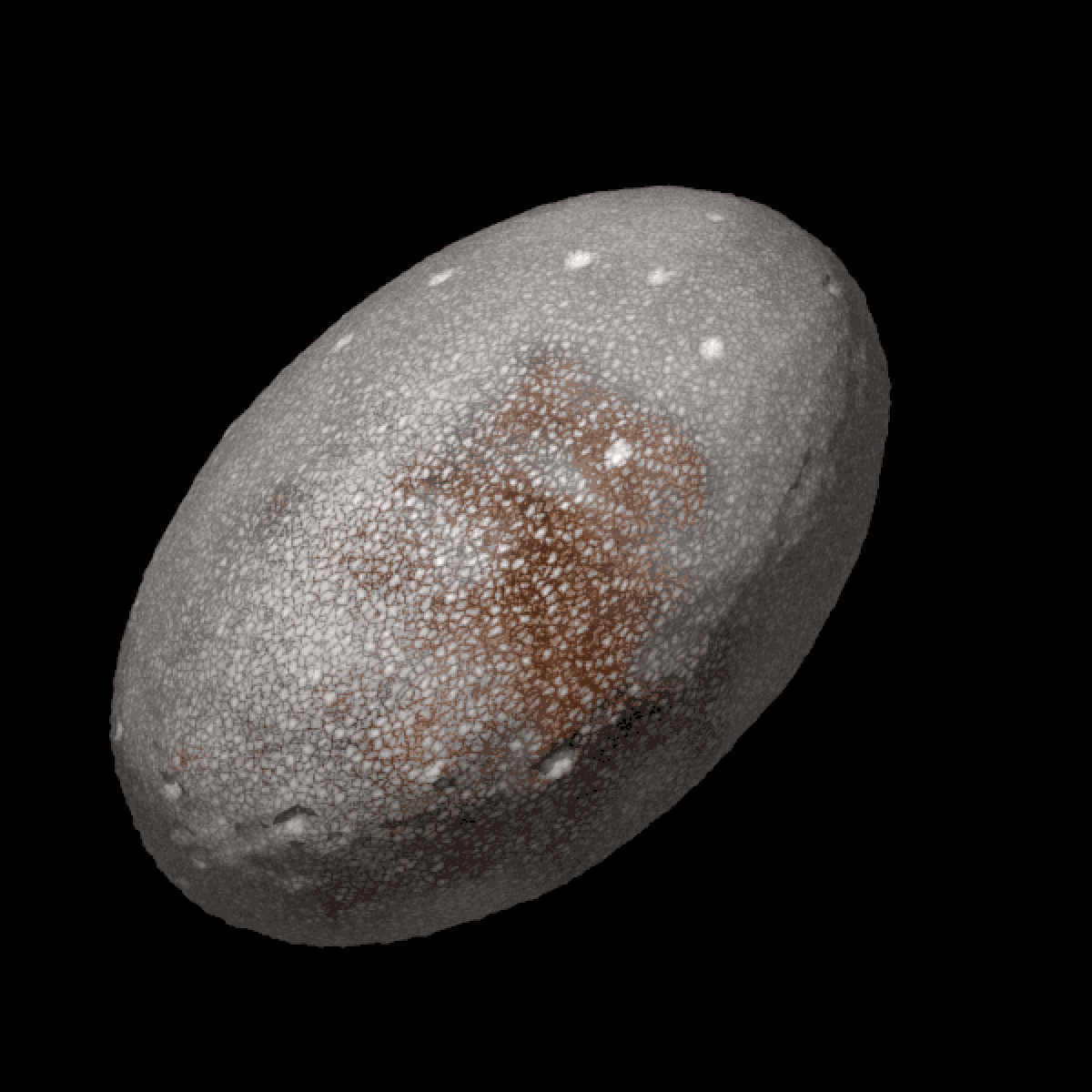}
  }  
  \subfloat[]{
    \includegraphics[height=2in]{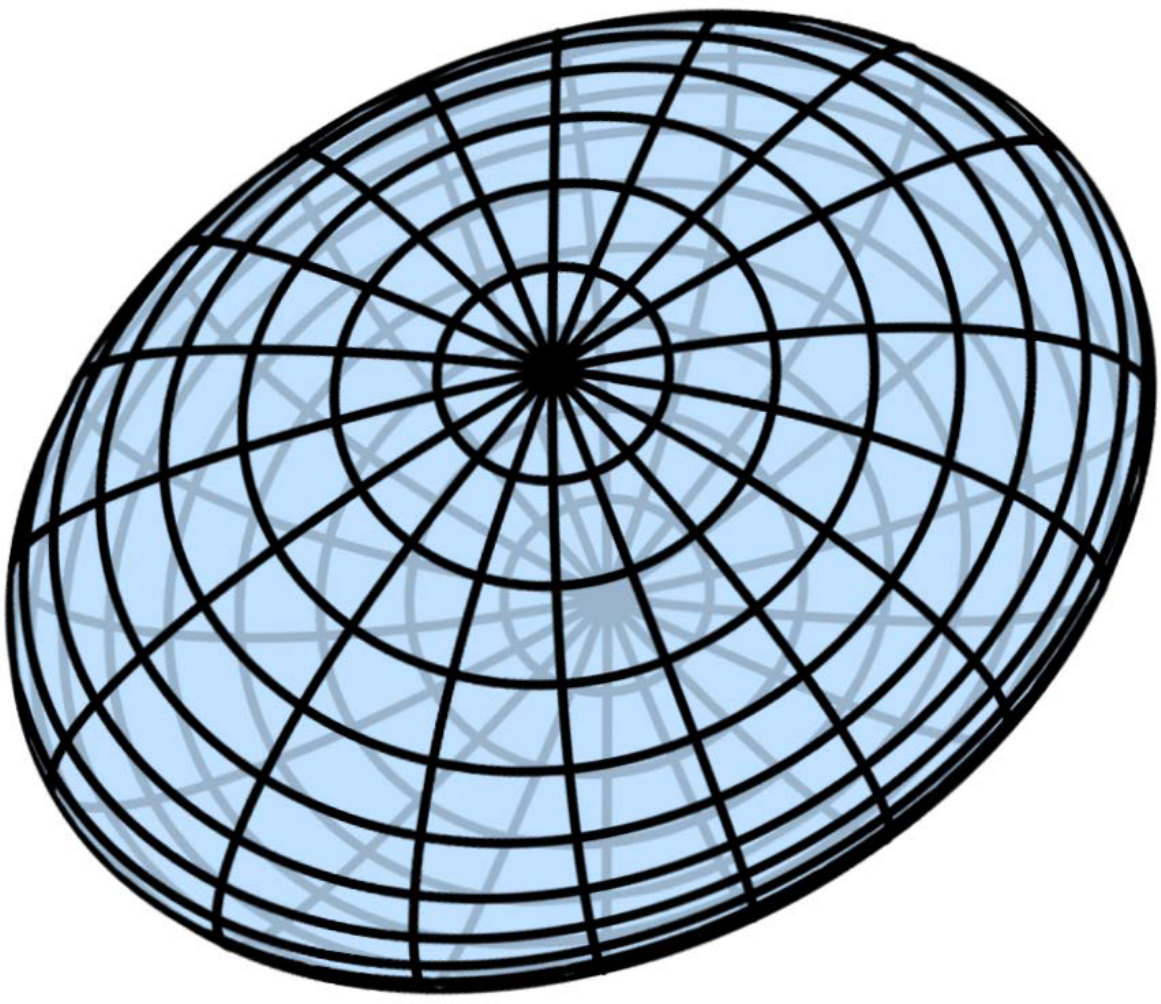}
  }\\
  \subfloat[]{
    \includegraphics[height=2.5in]{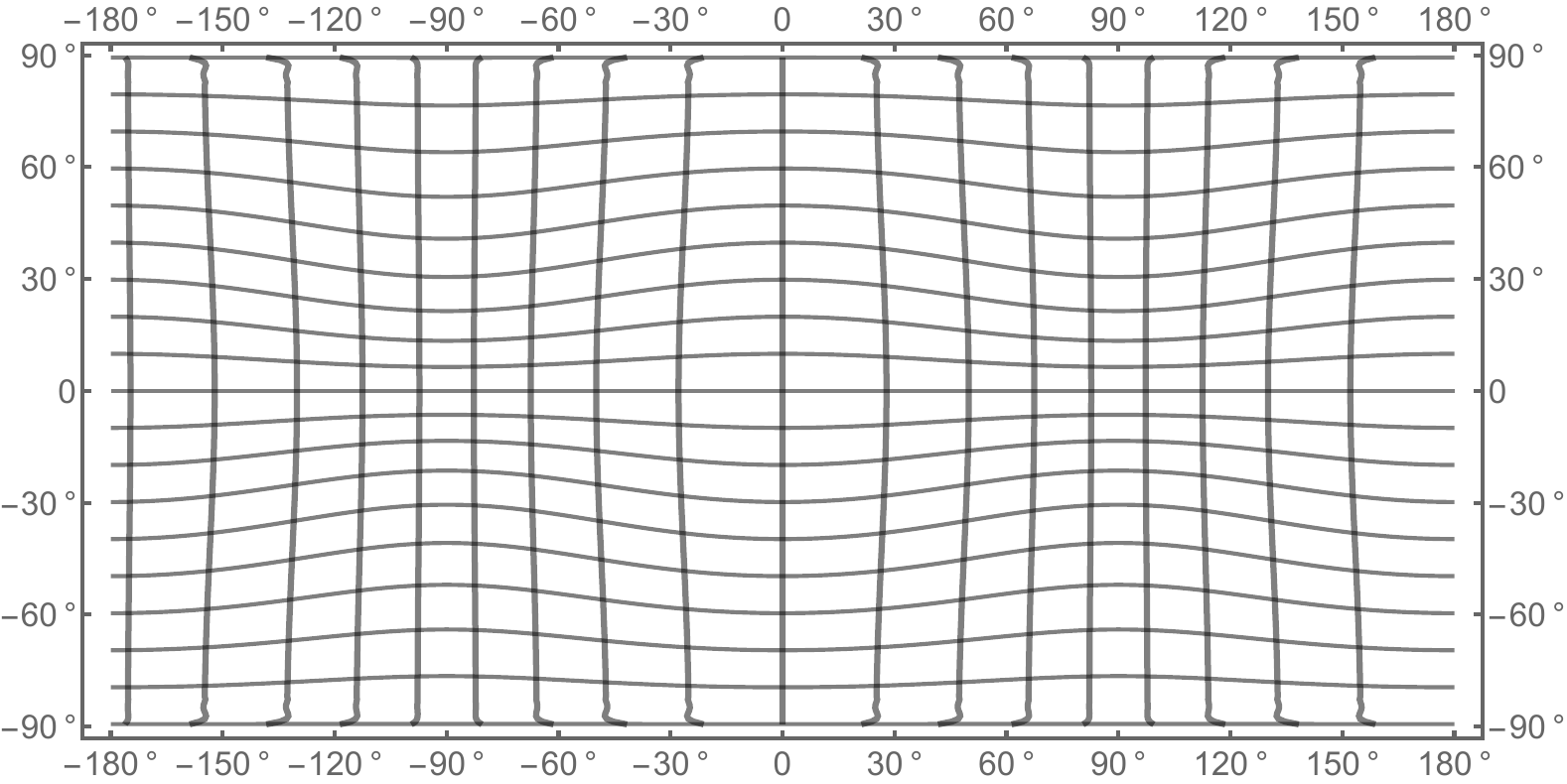}
  }
  \caption{Haumea. (a) Artist's impression, (b) latitudes and longitudes, and (c) world map.}
  \label{fig:appetizer}
\end{figure}

Obtaining a conformal mapping of an ellipsoid of a revolution onto a rectangular plane domain is a classical problem which arises from the construction of Mercator type map projections. Haumea is a dwarf planet located beyond Neptune's orbit discovered in the 2000's.
Due to its rapid rotation, its shape is a flattened ellipsoid with dimensions $1,960\times 1,518\times 996\,$km, which means that it is not
an ellipsoid of revolution. For a cartographer this means that in Mercator type map projections 
latitudes and longitudes do not intersect at right angles. This is illustrated in Figure~\ref{fig:appetizer}.
The algorithm used to compute the mapping is described in detail in the following sections.

\subsection{Organization}
The rest of the paper is organized as follows: After the introduction, preliminaries including the finite element method used and the
planar conjugate function method are covered. In short Section 3, the conjugate function method on surfaces is outlined.
In Section 4, two special surfaces, helicoid and catenoid are discussed, with a focus on their isothermal coordinates.
Numerical experiments, highlighting the efficacy of the method over a series of surfaces including a variety of
singularities, are the topic of the Section 5. Finally, conclusions are drawn in Section 6.

\section{Preliminaries}
\label{sec:preliminaries}
In this section, we review the preliminary material required for the development of the proposed method.
First, the planar conjugate function method is outlined.
Next, the $hp$-version of the finite element method as well as the two classes of a posteriori error
estimates are reviewed, including the a priori mesh refinement strategies
needed in the numerical experiments below.

\subsection{The Conjugate Function Method on Planar Domains}

We call a Jordan domain $\Omega \subset \mathbb{C}$ together with four distinct  oints $z_1,z_2,z_3,z_4\in \partial \Omega$ defining the positive orientation of $\partial \Omega$ a (generalized) {\it quadrilateral} and denote it by $Q=(\Omega; z_1,z_2,z_3,z_4)$. We denote the boundary segments connecting the pairs of points $(z_j,z_{j+1})$ for $j=1,2,3$, and $(z_4,z_1)$ for $j=4$, respectively, by  $\gamma_j$. 

It is well-known (see e.g. \cite{Ahlfors}) that there exists a unique number $h>0$ called the {\it conformal modulus} of $Q$, such that there exists a conformal mapping of the rectangle $R_h=[0,1]\times[0,h]\subset \C$ onto $\Omega$, with boundary points $z_1,z_2,z_3,z_4$ corresponding to the images of the points $0,1,1+ih,ih$, respectively. The number $h$ is unique and it determines the conformal equivalence class of $\Omega$ in the sense that there exists a conformal mapping between quadrilaterals (with boundary point correspondence) if and only if they have the same modulus. We denote the conformal modulus of a quadrilateral $Q$ by $M(Q)$.

The conformal modulus of a quadrilateral can also be determined as follows. Recall that there exists a (unique) harmonic solution $u$ to the following Dirichlet-Neumann mixed boundary value problem:
\begin{equation}
\label{dirneu}
\left\{
\begin{array}{rcl}
    \Delta u(z) = 0 & \text{for} & z \in \Omega, \\
    u(z) = 0 & \text{for} & z \in \gamma_2, \\
    u(z) = 1 & \text{for} & z \in \gamma_4, \\
    \partial u(z)/\partial n = 0 & \text{for} & z \in \gamma_1 \cup \gamma_3,
\end{array}
\right.
\end{equation}
where $n$ is the unit exterior boundary normal. Then the conformal modulus is connected to the above boundary value problem  by the identity (see e.g. Ahlfors \cite[Theorem~4.5]{Ahlfors} and Papamichael and Stylianopoulos \cite[Theorem~2.3.3]{PapamichaelStylianopoulos2010}):
\begin{equation}
M(Q)=\iint_{\Omega}|\nabla u|^2 \, dx \, dy.
\end{equation}

For a quadrilateral $Q=(\Omega; z_1,z_2,z_3,z_4)$ we call $\widetilde{Q}=(\Omega; z_2, z_3, z_4, z_1)$ its {\it conjugate quadrilateral} and the corresponding problem \eqref{dirneu} for quadrilateral $\widetilde{Q}$ the {\it conjugate Dirichlet--Neumann problem}. It is well-known that if $M(Q)=h>0$, then $M(\widetilde{Q})=1/h$, which leads to the following very useful {\it reciprocal identity}: For all quadrilaterals $Q$, we have
\begin{equation}
\label{resiproc}
M(Q) M(\widetilde{Q}) =1.
\end{equation}

Furthermore, we  may observe that the canonical conformal mapping of a quadrilateral $Q = (\Omega; z_1, z_2, z_3, z_4)$ onto the rectangle $R_h$ with vertices at $1+i h$, $i h$, $0$, and $1$, can be obtained by solving the corresponding Dirichlet--Neumann problem and its conjugate problem.

\begin{lemma}[\cite{hqr}]
Let $Q$ be a quadrilateral with modulus $h$, and suppose $u$ solves the Dirichlet--Neumann problem \eqref{dirneu}. If $v$ is a harmonic function conjugate to $u$, satisfying $v(\operatorname{Re} z_3, \operatorname{Im} z_3)=0$, and $\tilde{u}$ represents the harmonic solution for the conjugate quadrilateral $\widetilde{Q}$, then $v = h \tilde{u}$.
\end{lemma}

This lemma gives a method for computing numerical conformal mappings between plane domains, which we call the {\it conjugate function method}. The accuracy of this method in the planar case is discussed in \cite[pg. 348]{hqr}.

\subsection{High-Order Finite Element Method}\label{sec:hpfem}

High-order finite element methods have the capability for exponential convergence
provided the discretisation is constructed properly in both domain ($h$-version)
and in local polynomial order ($p$-version). In this paper, the combined $hp$-FEM
is used to discretise the Laplace-Beltrami operator in the parameter space
of the surfaces considered.

In all cases, it is implicitly assumed that the exact parameterisation of
the boundaries on the parameter space is known. This allows us to benefit
from efficient handling of large elements within the $p$-version
without significant loss of accuracy. It also means that the number of 
elements can be kept relatively low.

Let us consider the Dirichlet-Neumann problem \eqref{dirneu} and its weak solution $u_0$.
The following theorem due to Babu{\v{s}}ka and Guo \cite{BaGuo1,BaGuo2},
sets the limit to the rate of convergence. Notice that the construction of
the appropriate spaces is technical, but can be extended to parameterised surfaces. 
For rigorous treatment of the theory involved
see Schwab \cite{schwab} and references therein. 

\begin{theorem} \label{propermesh}
Let $\Omega \subset \mathbb{R}^2$ be a polygon, $v$ the FEM-solution of \eqref{dirneu}, and
let the weak solution $u_0$ be in a suitable countably normed space where
the derivatives of arbitrarily high order are controlled.
Then
\[
\inf_v \|u_0 - v\|_{H^1(\Omega)} \leq C\,\exp(-b \sqrt[3]{N}),
\]
where $C$ and $b$ are independent of $N$, the number of degrees of freedom. 
Here $v$ is computed on a proper geometric mesh, where the order of an individual
element is set to be its element graph distance to the nearest singularity.
(The result also holds for meshes with constant polynomial degree.)
\end{theorem}

\subsubsection{Error Estimation}
We employ two a posteriori error estimation methods. 
The first one is the direct application of the
reciprocal equation \eqref{resiproc} as an error estimate. We simply monitor the
quantity of interest
\begin{equation}
 \operatorname{reci}(Q) = | M(Q) M(\widetilde{Q}) - 1 |,
\end{equation}
which is defined for every quadrilateral $Q$ and will be shown to be applicable also on surfaces.
In fact, existence of this type of relations can even be used as necessary and sufficient
conditions for the existence of quasiconformal parameterization of a surface in a very general setting
\cite[Theorem 1.4]{Rajala} of which our results are a special case.

Similarly, the auxiliary subspace error estimation can be defined analogously
to Laplace-Beltrami. Let $\cT$ be some $hp$-discretisation on the
computational domain $\Omega$. Assuming that the exact
solution $u \in H_0^1(\Omega)$, defined on $\cT$, has finite energy, the
approximation problem is as follows: Find $\hat{u} \in V$ such that
\begin{equation}\label{eq:approximation}
  a(\hat{u},v)=l(v)\ (= a(u,v))\quad (\forall v \in V),
\end{equation}where
$a(\,\cdot\,,\,\cdot\,)$ and $l(\,\cdot\,)$, are the bilinear form
and the load potential, respectively. Additional degrees of
freedom can be introduced by enriching the space $V$. This is
accomplished via the introduction of an auxiliary subspace or ``error
space'' $W \subset H_0^1(\Omega)$ such that $V \cap W = \{0\}$. 
The error problem becomes thus: Find $\varepsilon \in W$ such that
\begin{equation}\label{eq:error}
  a(\varepsilon,v)=l(v)- a(\hat{u},v) (= a(u-\hat{u},v))\quad (\forall v \in W).
\end{equation}
This can be interpreted as a projection of the residual to the auxiliary space.

The main error theorem on auxiliary subspace error estimators for 
standard diffusion problems is Theorem~\ref{KeyErrorThm}.
It should be mentioned that even though there exists compelling
numerical evidence that the constant $K$ is in fact independent of $p$,
no rigorous proofs exist to support this observation.

\begin{theorem}[\cite{hno}]\label{KeyErrorThm}
There is a constant $K$ depending only on the dimension $d$,
polynomial degree $p$, continuity and coercivity constants $C$ and
$c$, and the shape-regularity of the triangulation $\cT$ such that
\begin{align*}
\frac{c}{C}\,\|\eE\|_{H^1(\Omega)}\leq\|u-\hat{u}\|_{H^1(\Omega)}\leq K
\left(\|\eE\|_{H^1(\Omega)}+\osc(R,r,\cT)\right),
\end{align*}
where the residual oscillation depends on the volumetric and face
residuals $R$ and $r$, and the triangulation $\cT$.
\end{theorem}

\subsubsection{A Priori Refinement Strategies}\label{sec:refine}
\begin{figure}
  \centering
  \subfloat[]{\includegraphics[height=1.5in]{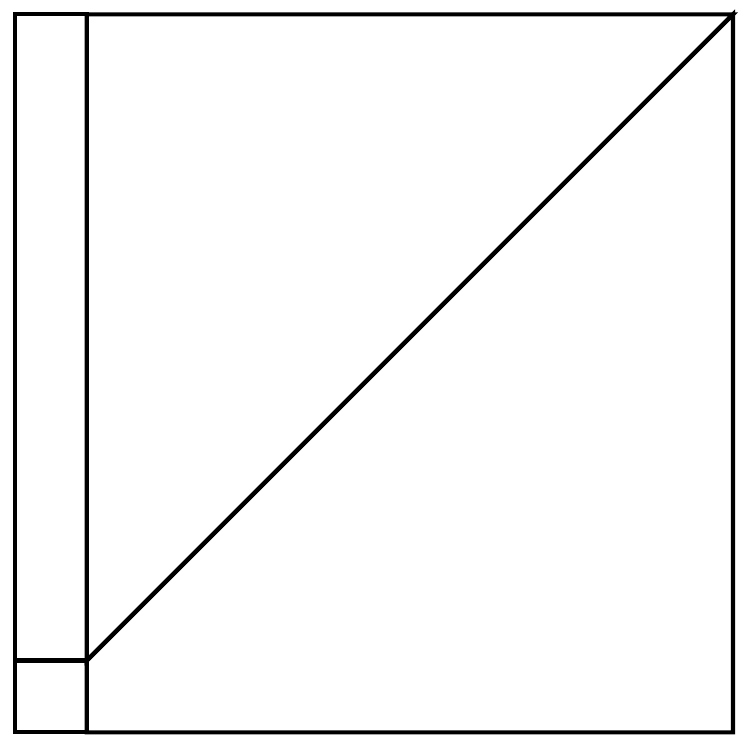}}\quad
  \subfloat[]{\includegraphics[height=1.5in]{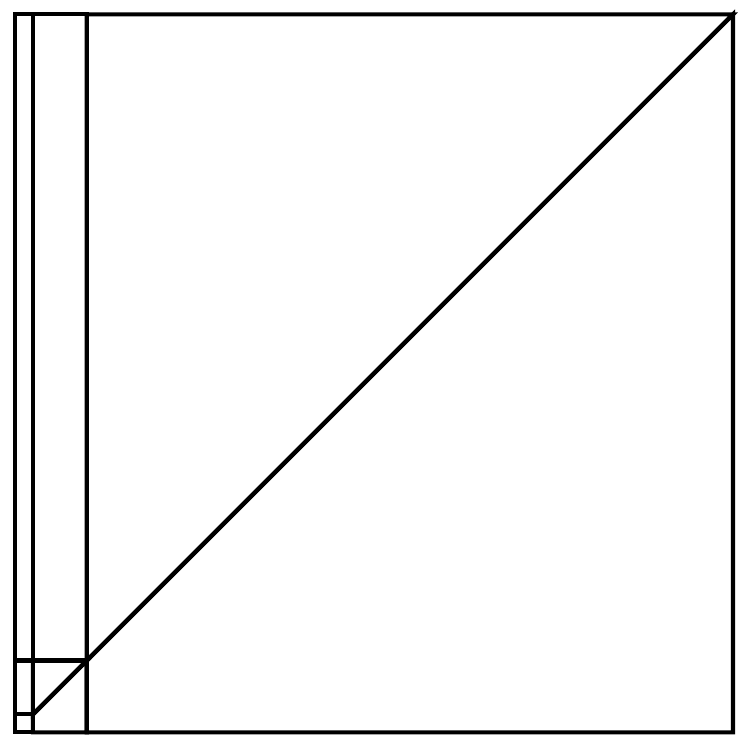}}
  \caption{Rule based mesh refinement. Refinement to an edge coupled with corner refinement.
(a) and (b): The first two steps starting from a regular grid. Notice that the process generates both triangles and quadrilaterals.}\label{fig:refinement}
\end{figure}

In many problems exponential convergence of the $p$- and $hp$-FEM cannot be realised unless
conformal meshes with geometric grading are available.
Our approach is based on generating a priori optimally refined meshes
using rule based algorithms \cite{Hakula2013}. 
One of the interesting requirements here is that 
we have to find a way to simultaneously refine the mesh geometrically over the
whole edge or segment of the boundary coupled with standard corner refinement.
The replacement rules are illustrated in Figure~\ref{fig:refinement}.

Since the geometric grading strategy relies on a priori information
it does not lend itself well for adaptivity. Modifying the mesh in a posteriori sense
requires changing the a priori strategy, in other words, remeshing the whole domain.
Hence, we only compute error estimates as defined above
and do not propose any algorithmic approach for adaptivity.

\section{Conjugate Function Method on Surfaces}
\label{sec:cfm}

In this section, we develop a version of the conjugate function method introduced in \cite{hqr,hqr2} in the case of domains on Riemannian surfaces. For the case of simply connected surfaces, the existence of the harmonic conjugate is well-known \cite[Corollary 9-2, pg. 135]{Ahlfors}.
Our discussion covers the variational formulation used in the
finite element method.

\subsection{Conformal Modulus and Generalized Quadrilaterals of Surfaces}
Suppose $\Omega\subset \Sf$ is a simply connected domain on a Riemannian surface $\Sf$ so that its boundary of $\Omega$ is a Jordan curve.
Let $z_1,z_2,z_3,z_4\in \partial \Omega$ be four distinct points on the boundary of $\Omega$ defining the positive orientation of  $\partial \Omega$. We again call $Q=(\Omega; z_1,z_2,z_3,z_4)$ a quadrilateral on $\Sf$, and define its conformal modulus as the number $h>0$ such that there exists a conformal mapping of the rectangle $R_h=[0,1]\times[0,h]\subset \C$ onto $\Omega$. with boundary points $z_1,z_2,z_3,z_4$ corresponding to the images of the points $0,1,1+ih,ih$, respectively. The number $h$ is unique and it determines the conformal equivalence class of $\Omega$.

Denote by $\gamma_j$ the part of the boundary curve $\partial \Omega$ connecting $z_j,z_{j+1}$ for $=1,2,3$ and $z_4,z_1$ for $j=4$. As in the planar case (see \cite[2.1]{hrv}), numerical computation of the conformal modulus can be based on the following mixed boundary value problem:
\[
\Delta_{\Sf}u=0\text{ on }\Omega,
\]
with Dirichlet boundary conditions $0$ on $\gamma_2$ and $1$ on $\gamma_4$, respectively, and Neumann boundary condition $\partial u/\partial n=0$ on $\gamma_1\cup\gamma_3$, where $n$ is the exterior unit normal at a boundary point taken at the tangent plane of $\Sf$. Numerical examples on computation of conformal moduli of plane quadrilaterals are given in \cite{hrv}. 

\subsection{Laplace-Beltrami}\label{sec:lb}
Our task is to define the Laplacian on some surface $S$, that is, we want to define the operator
$\Delta_S$ in the form which is suitable for finite element implementation.
In our setting, the surface is always assumed to be
given in some parameterised form. Let $\mathbf{x}_S : \bar{S} \to S$ be a parameterisation
of a surface $S$. The goal is to treat $\bar{S} \subset \mathbb{R}^2$ as the reference domain on
which the finite elements are defined. Let $J_\mathbf{x}$ be the Jacobian of the mapping, and
hence $G_S = J_\mathbf{x}^T J_\mathbf{x}$ is the first fundamental form.

The tangential gradient of some function $v : S\to \mathbb{R}$ is
\begin{equation}
  (\nabla_S v)\circ \mathbf{x}_S := J_\mathbf{x} G_S^{-1} \nabla(v \circ \mathbf{x}_S),
\end{equation}
and immediately, using the same notation, the $\Delta_S$ can be written as
\begin{equation}
  \Delta_S := \nabla_S \cdot \nabla_S.
\end{equation}

Equipped with this operator, the variational formulation
\begin{equation}
\int_S \nabla_{S}\psi\cdot\nabla_{S} v\,dx,\ \mbox{ for all }v\in H^1(S)
\end{equation}
on an image $K$ of a given element $\bar{K}$ in a discretisation of $\bar{S}$ becomes
\begin{align}\label{eq:femvar}
\int_K \nabla_{K}\psi\cdot\nabla_{K} v\,dx &= 
\int_{\bar{K}} \nabla(\psi \circ \mathbf{x}_K)^T  G_K^{-T} J_{\mathbf{x}}^T J_\mathbf{x} G_K^{-1} \nabla(v \circ \mathbf{x}_K) \sqrt{\operatorname{det}(G_k)}\,d\bar{x} \\
&= \int_{\bar{K}} \nabla(\psi \circ \mathbf{x}_K)^T  G_K^{-T} G_S G_K^{-1} \nabla(v \circ \mathbf{x}_K) \sqrt{\operatorname{det}(G_k)}\,d\bar{x}.
\end{align}

The integrals are evaluated on standard 2D mapped Gaussian quadratures. In fact, it is the first
equality \eqref{eq:femvar} which is compatible with our implementation of the method.
Since the problem has been transformed to a standard 2D planar problem with variable
coefficients, there is no need for additional arguments on the convergence of the method.

\section{Helicoid and Catenoid: Isothermal Coordinates}
\begin{figure}
  \centering
  \subfloat{
    \includegraphics[height=2.5in]{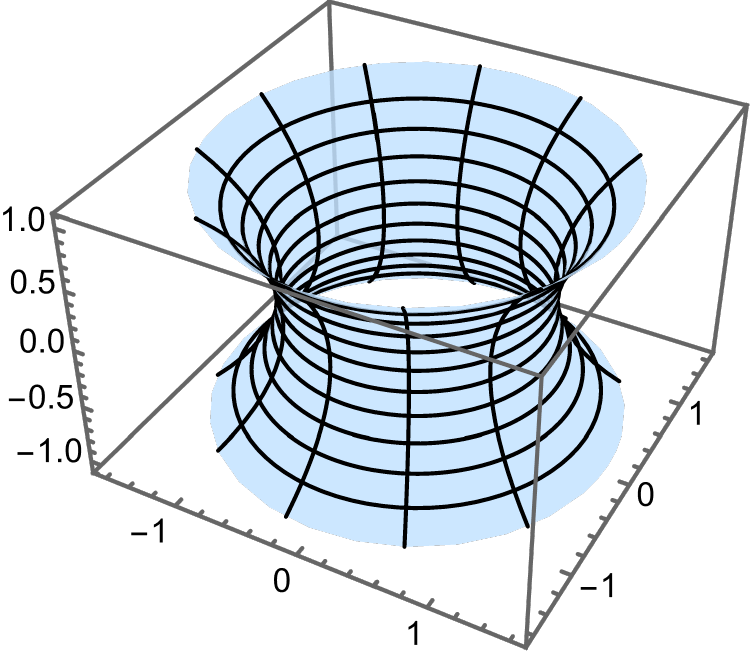}
  }\quad
  \subfloat{
    \includegraphics[height=2.5in]{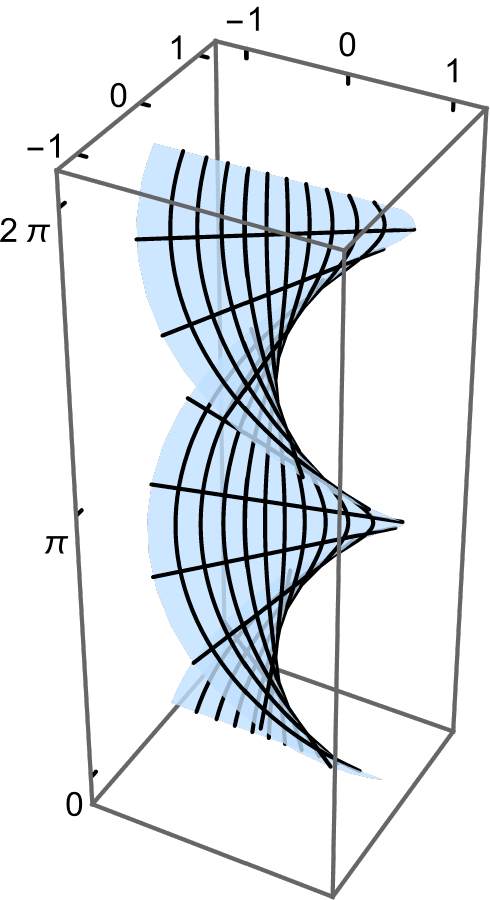}
  }
  \caption{Catenoid and Helicoid: Conformal maps on the surfaces.}
  \label{fig:isothermal}
\end{figure}

Isothermal coordinates on a Riemannian manifold are local coordinates 
where the metric is conformal to the Euclidean metric. 
In the context of our study, this leads to two immediate observations:
First, the conformal moduli can be evaluated exactly, if the parameter space is one of the known
conformal domains; second, in the variational formulation \eqref{eq:femvar} the surface
contribution to the mapping
$G_K^{-T} G_S G_K^{-1} \sqrt{\operatorname{det}(G_K)}$
is an identity mapping if the parameterisation is isothermal, thus giving us a simple verification test.

\subsection{Catenoid}
Consider the graph $x=\cosh(z)$ on the $xz$-plane, which is called {\it catenary}. The surface of revolution generated by rotation of this graph around the $z$-axis is called {\it catenoid}, and its isothermal coordinates are given by
\begin{equation}\label{eq:cate-isotherm2}
\vr_C(u,v)=\big(x(u)\cos(v),x(u)\sin(v),z(u)\big),
\end{equation}
where $(x(u),z(u))=(\cosh(u),u)$. 
Actually, in this case the first fundamental form is given by $\cosh^2(u)(du^2+dv^2)$, 
in other words, the metric is conformal to the Euclidean metric.
The mapping information used in the variational formulation is
\[
J_\mathbf{x} =\left(
\begin{array}{cc}
 \sinh (u) \cos (v) & -\cosh (u) \sin (v) \\
 \sinh (u) \sin (v) & \cosh (u) \cos (v) \\
 1 & 0 \\
\end{array}
\right), \quad
G_K^{-1}=\left(
\begin{array}{cc}
 \text{sech}^2(u) & 0 \\
 0 & \text{sech}^2(u) \\
\end{array}
\right),
\]
\[
\sqrt{\operatorname{det}(G_K)} =\cosh ^2(u).
\]
Through a straightforward computation one finds that, indeed, 
\[
G_K^{-T} G_S G_K^{-1} \sqrt{\operatorname{det}(G_K)} = I.
\]
\subsection{Helicoid}
Recall that a parameterisation of a {\it helicoid} surface is given by
\begin{equation}\label{heli-non-isotherm}
\vr_{H_1}(u,v)=\big(u \cos(v), u\sin(v),v\big).
\end{equation}
This parameterisation is not isothermal. On the other hand, an isothermal parameterisation of helicoid 
exists (see \cite{Oprea}), and can be written as 
\begin{equation}\label{heli-isotherm}
\vr_{H_2}(u,v)=\big(\sinh(u) \sin(v), -\sinh(u)\cos(v),v\big).
\end{equation}

As for the catenoid above, for $\vr_{H_1}(u,v)$, one has
\[
J_\mathbf{x} =\left(
\begin{array}{cc}
 \cos (v) & -u  \sin (v) \\
  \sin (v) & u  \cos (v) \\
 0 & 1 \\
\end{array}
\right), \quad
G_K^{-1}=\left(
\begin{array}{cc}
 1 & 0 \\
 0 & \frac{1}{u^2 + 1} \\
\end{array}
\right),
\]
\[
\sqrt{\operatorname{det}(G_K)} =\sqrt{u^2+1},
\]
resulting in
\begin{equation}\label{eq:nonisothermal}
G_K^{-T} G_S G_K^{-1} \sqrt{\operatorname{det}(G_K)} = \left(
\begin{array}{cc}
 \sqrt{u^2+1} & 0 \\
 0 & \frac{1}{\sqrt{u^2+1}} \\
\end{array}
\right) \neq I.
\end{equation}

For $\vr_{H_2}(u,v)$, one has
\[
J_\mathbf{x} = \left(
\begin{array}{cc}
 \cosh (u) (-\cos (v)) & \sinh (u) \sin (v) \\
 \cosh (u) \sin (v) & \sinh (u) \cos (v) \\
 0 & 1 \\
\end{array}
\right), \quad
G_K^{-1}=\left(
\begin{array}{cc}
 \text{sech}^2(u) & 0 \\
 0 & \text{sech}^2(u) \\
\end{array}
\right),
\]
\[
\sqrt{\operatorname{det}(G_K)} =\cosh ^2(u),
\]
and the mapping in the variational formulation reduces to identity.

\subsection{Computation of Moduli: Exact Values and Error Estimates}
For both helicoid and catenoid the parameter domain $\Omega = [-1,1]\times [0,2\pi]$.
The exact values of the conformal moduli for both quadrilaterals with natural corner points are
\begin{equation}
  M(Q) = 1/\pi, \quad M(\widetilde{Q}) = \pi.
\end{equation}

\begin{figure}
  \centering
  \subfloat[{$(u,v) \in [-1,1]\times [0,2\pi]$}]{\label{fig:helicoidA}
    \includegraphics[height=2in]{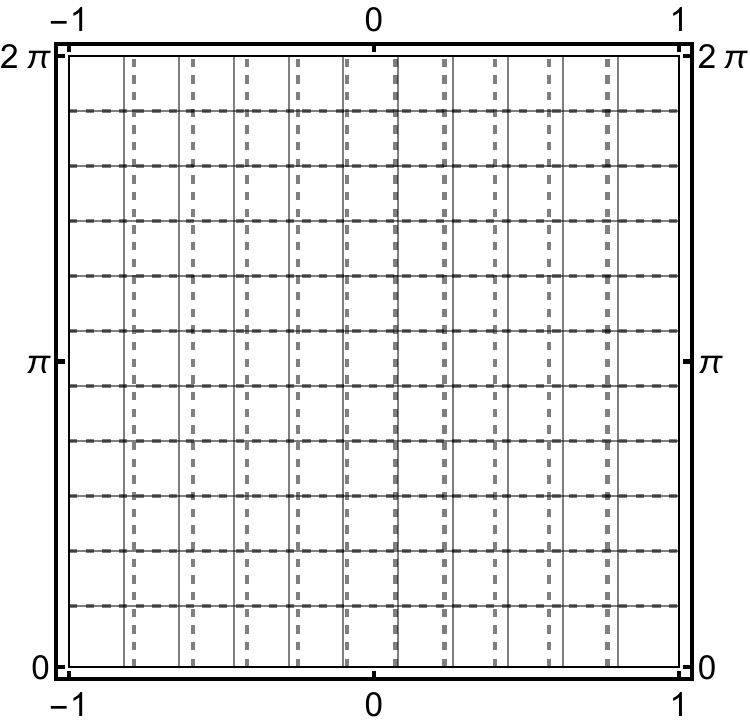}
  }
  \subfloat[Error estimates vs. $p$.]{\label{fig:helicoidB}
    \includegraphics[width=0.55\textwidth]{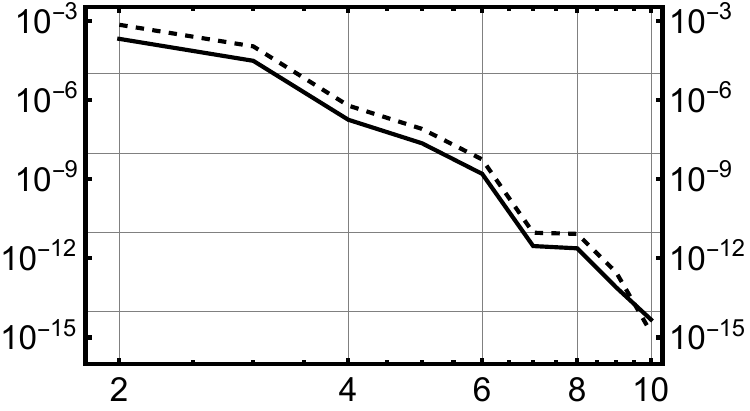}
  }
  \caption{Helicoid: (a) Maps on the parameter space. Isothermal coordinates: Solid lines, General: Dashed lines. (b) Convergence in the non isothermal parameterisation, error estimates vs. $p$. Reciprocal error estimate: Solid line,
Auxiliary subspace error estimate: Dashed line. 
}
  \label{fig:helicoidcoordinates}
\end{figure}

The effect of the parameterisation is illustrated in Figure~\ref{fig:helicoidA}.
The conformal mapping is perfectly aligned with the coordinate axes if the parameterisation is
isothermal, whereas in the general case of $\vr_{H_1}(u,v)$ (dashed lines) there is some
deviation in the $u$-direction. Notice, that in the $v$-direction there is perfect alignment,
as one would expect since in \eqref{eq:nonisothermal} there is no dependence on $v$.

For isothermal parameterisations linear finite elements are exact. For $\vr_{H_1}(u,v)$
in the $u$-direction this cannot be the case. In Figure~\ref{fig:helicoidB} convergence graphs
for both reciprocal and auxiliary space error estimators are shown. The convergence is exponential
in $p$, and interestingly exhibits mild staircasing.


\section{Numerical Experiments}
\label{sec:numex}

The numerical experiments are, broadly speaking, divided into three sections.
First we continue the discussion on cartography and consider the planet Earth, since remarkably
the conformal mapping is known and even pointwise convergence can be estimated.
Second, we study problems on parts of a sphere with different characteristic features, such as singularities
induced by the boundary conditions, cusps, or higher connectivity of the domain.
In two specific problems on the sphere, the exact values of the conformal moduli are known, and hence
the convergence rates of the method in various natural norms can be shown to be exponential
provided the meshes are properly refined.
As our last example, a highly complicated non-smooth surface is considered.
This section is concluded with remarks on computational complexity.
\subsection{Planet Earth (Ellipsoid)}\label{sec:earth}
Using values for the equatorial and polar radii of the planet Earth (as specified by
World Geodetic System, 1984), $R_1 = 6378.1370$ and $R_2 = 6356.7523$, respectively, 
we arrive at the following normalised parametrisation, where the parameters are
the longitude $\lambda$ and latitude $\phi$
\begin{equation*}\label{eq:ellipsoid}
\left\{\begin{array}{rcl}
x &=&a \cos (\lambda ) \sin (\phi ), \\
y &=&b \sin (\lambda ) \sin (\phi ), \\ 
z &=&c \cos (\phi ),
\end{array}\right.
\end{equation*}
with $a = b = R_1/R_2$, and $c=1$. For the reference case, that is, the globe with $a = b = c = 1$, the Mercator 
projection can be computed in closed form, where a point $(\xi,\eta)$ on a rectangular map is given as a function
of its longitude and latitude \cite[pg. 24]{Vermeer-Rasila}:
\begin{equation}\label{eq:mercator}
  (\xi,\eta) = \Big(\lambda, \ln \tan\Big(\frac{\pi}{4} + \frac{\phi}{2}\Big)\Big).
\end{equation}
As is well-known, the projection is singular at the poles, hence we remove small polar caps of radius $\epsilon > 0$,
and the parameter domain is thus $\lambda \in [0,2\pi]$, $\phi \in [\epsilon,\pi-\epsilon]$. We choose $\epsilon = 1/100$.

\begin{figure}
  \centering
  \subfloat[]{
    \includegraphics[width=0.45\textwidth]{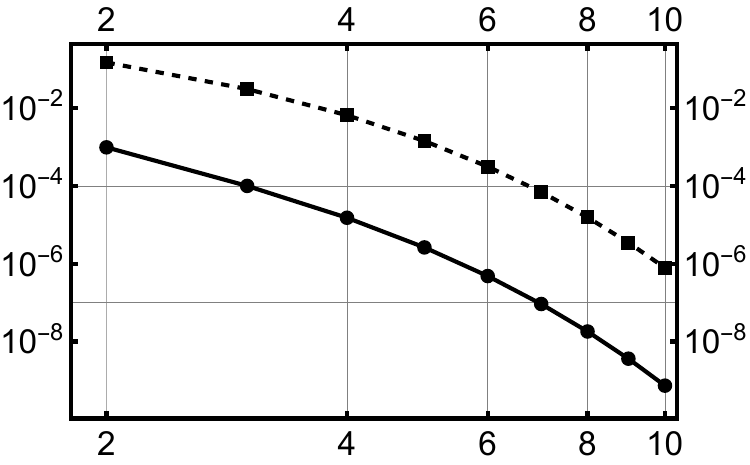}
  }
  \subfloat[]{
    \includegraphics[width=0.45\textwidth]{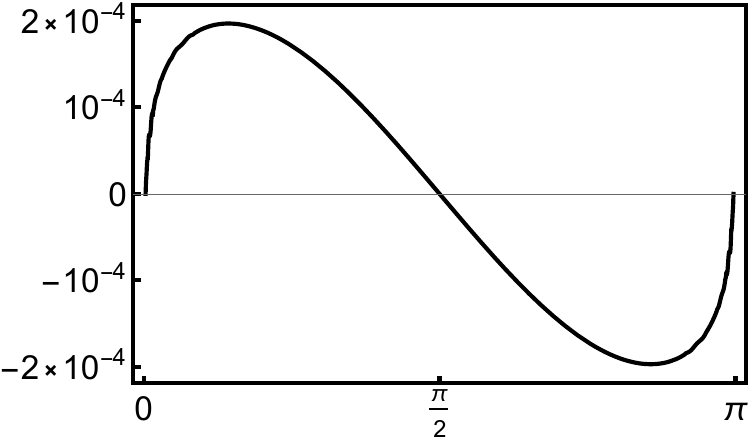}
  }
  \caption{Mercator projection. 
  (a) Convergence measured in $L^2$-norm (solid line) and $H^1$-seminorm (dashed line) vs polynomial order (loglog-plot). 
  (b) Effect of flattening. Difference between the normalised Mercator projection and the observed one for the
Earth ellipsoid along $\lambda = 0$.}
  \label{fig:mercator}
\end{figure}

For both cases the proposed method is accurate as in the cases above as measured in the reciprocal error. Since the
exact mapping is known, we can measure error also over the parameter domain. 
We normalise the projection \eqref{eq:mercator} to have values in the interval $[0,1]$ along $\lambda = 0$
and compute $L^2$-norm and $H^1$-seminorm errors using the $hp$-FEM solution.
With a boundary layer fitted mesh at $p=10$ the norms are $\sim 10^{-9}$ and $\sim 10^{-6}$
for $L^2$-norm and $H^1$-seminorm, respectively. The convergence graphs are shown in Figure~\ref{fig:mercator} 
as well as the effect of the Earth's flattening
on the mapping, which is illustrated with the difference between the observed and the reference solutions.
The effect is small, but visible.


\begin{table}
\centering
\caption{Surfaces studied in the numerical experiments. Surface parameterisation is: $\left(\sin (u) \cos (v),\sin (u) \sin (v),\cos (u)\right)$, i.e., a sphere.
For the hyperbolic quadrilateral and multiply connected examples the parameter domains are nontrivial.}
\label{tbl:experiments}

\begin{tabular}{ll}
Name & Parameter Domain\\ \hline
Hemisphere               & $u\in [0,\pi/2], v\in [0,2\pi]$\\
Quarter Sphere           & $u\in [0,\pi/2], v\in [0,\pi]$\\
Hyperbolic Quadrilateral & $\text{HypQuad}((3\pi/16,3\pi/16),\frac{\pi }{8 \sqrt{2}},\pi/4)$\\
Two Holes                & $\mathbb{D}((\frac 12, \frac 12),1)\setminus (\mathbb{D}((\frac {1} {4},\frac 14),\frac 14)
\cup \mathbb{D}((\frac 34,\frac 34),\frac 14))$ \\
\end{tabular}
\end{table}

\subsection{Schwarzian Hemisphere}
\begin{figure}
  \centering
  \subfloat{
    \includegraphics[height=2.5in]{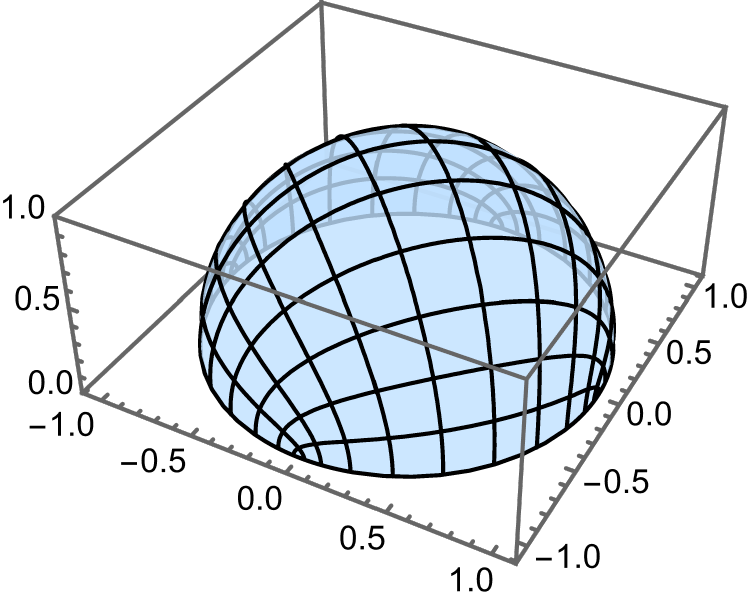}
  }
  \caption{Schwarzian Hemisphere. Mapping on the surface.}
  \label{fig:schwarzgrid}
\end{figure}

Our first example with singularities is a surface variant of the famous example by Schwarz.
The rim of a hemisphere is divided into four sections of equal length, in other words, the
four corners define a quadrilateral. The resulting map is shown in Figure~\ref{fig:schwarzgrid}.

The four points $z_i$ on the boundary are
\[
z_1 = (1,0,0),z_2=(0,1,0),z_3=(-1,0,0),z_4=(0,-1,0).
\]
Due to symmetry it follows immediately that $M(Q) = M(\widetilde{Q}) = 1$.

\begin{figure}
  \centering
  \subfloat[]{\label{fig:schwarzfemA}
    \includegraphics[height=2in]{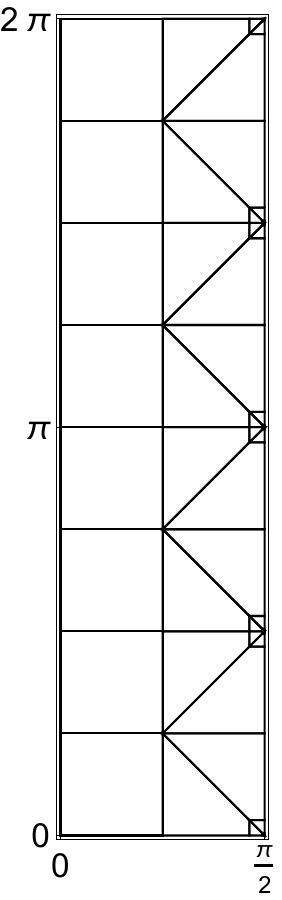}
  }\quad
  \subfloat[]{\label{fig:schwarzfemB}
    \includegraphics[height=1.975in]{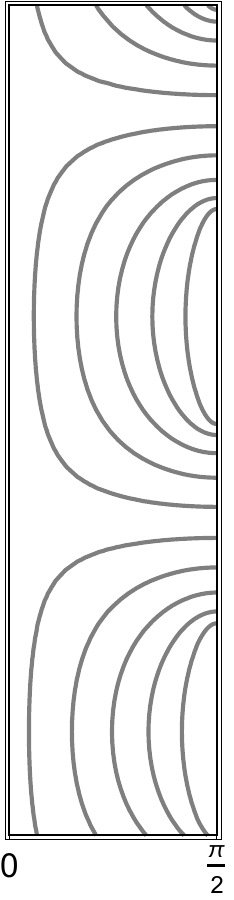}
  }\quad
  \subfloat[]{\label{fig:schwarzfemC}
    \includegraphics[height=1.975in]{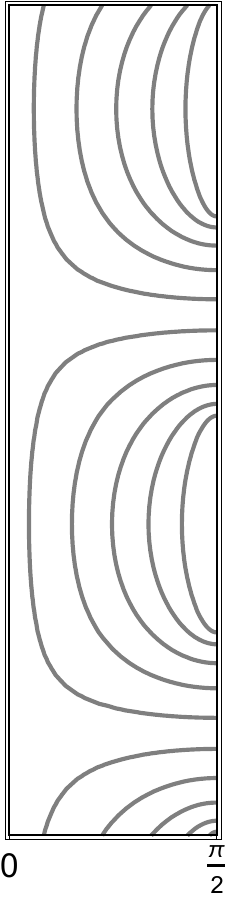}
  }\quad
  \subfloat[]{\label{fig:schwarzfemD}
    \includegraphics[height=2in]{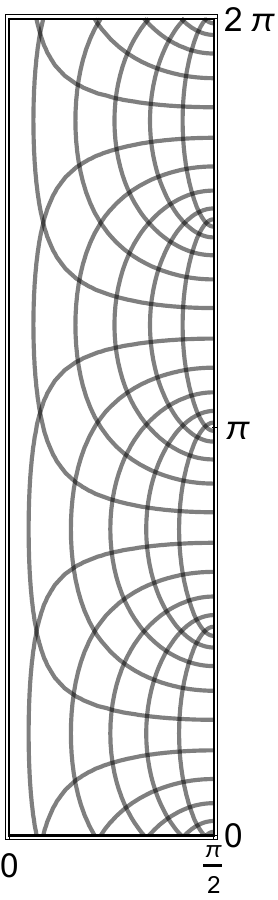}
  }
  \caption{Schwarzian Hemisphere. Construction on the parameter plane. (a) Mesh with geometric grading to all singularities.
  (b), (c) Contour lines of the solutions of the original and conjugate problems. (d) Projection of the conformal map onto
the parameter domain.}
  \label{fig:schwarzfem}
\end{figure}

The construction on the parameter plane is illustrated in Figure~\ref{fig:schwarzfem}.
The locations of the symmetries are the same for both the original and 
conjugate problems and therefore only one mesh is needed (see Figure~\ref{fig:schwarzfemA}).
The contour lines of the respective solutions are shown in Figures~\ref{fig:schwarzfemB}-~\ref{fig:schwarzfemC},
and brought together in Figure~\ref{fig:schwarzfemD}. It is only after the map is lifted onto the actual surface
when the contour lines become locally orthogonal.

The observed exponential convergence is shown in Figure~\ref{fig:sphereconvergenceA}.
Three types of errors are shown: The reciprocal error, the exact error for $M(Q)$, and
its estimated error.
Since the reciprocal error measures the error of a product, it should be an upper bound for both factors in the product.
This is indeed the case, also the estimated error is very accurate and remains asymptotically consistent.

\subsection{Quarter Sphere}
\begin{figure}
  \centering
  \subfloat{
    \includegraphics[height=2.5in]{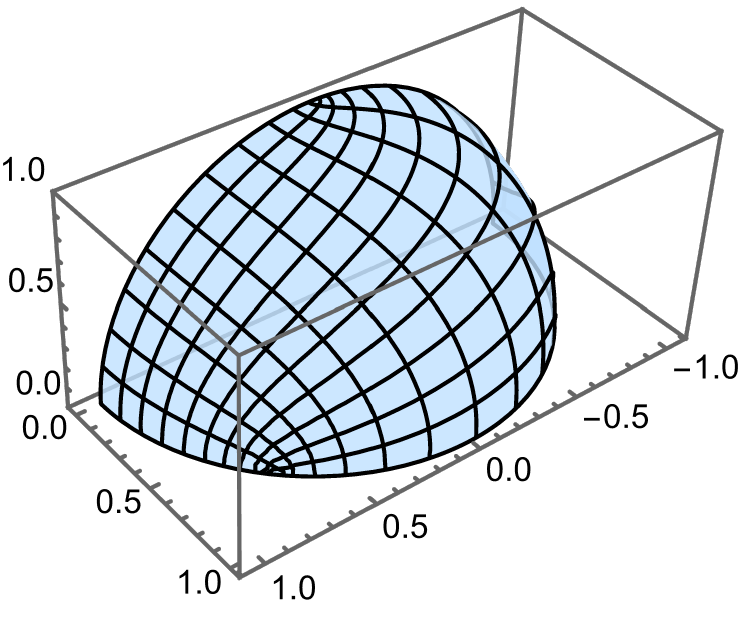}
  }
  \caption{Quarter Sphere. Mapping on the surface}
  \label{fig:qsgrid}
\end{figure}

\begin{figure}
  \centering
  \subfloat[]{\label{fig:qsfemA}
    \includegraphics[height=2in]{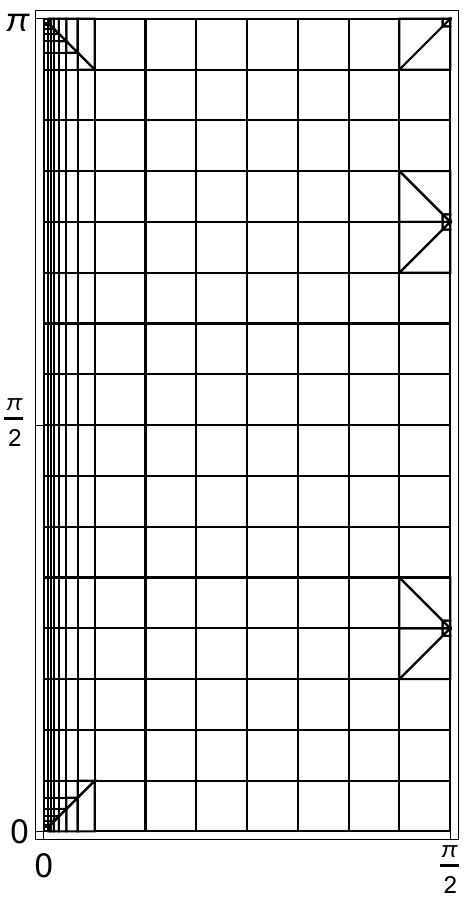}
  }
  \subfloat[]{\label{fig:qsfemB}
    \includegraphics[height=1.975in]{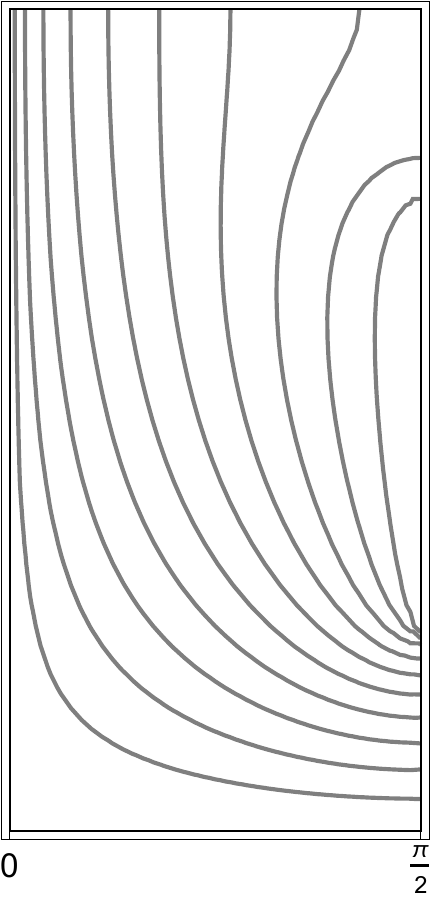}
  }
  \subfloat[]{\label{fig:qsfemC}
    \includegraphics[height=1.975in]{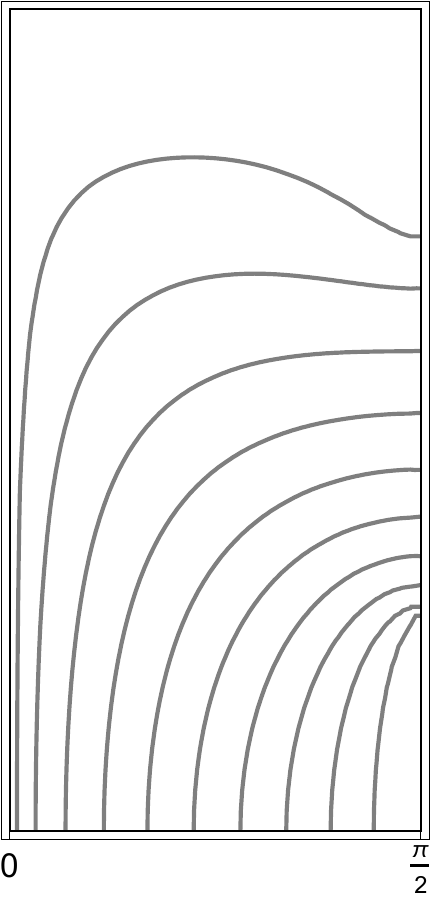}
  }
  \subfloat[]{\label{fig:qsfemD}
    \includegraphics[height=2in]{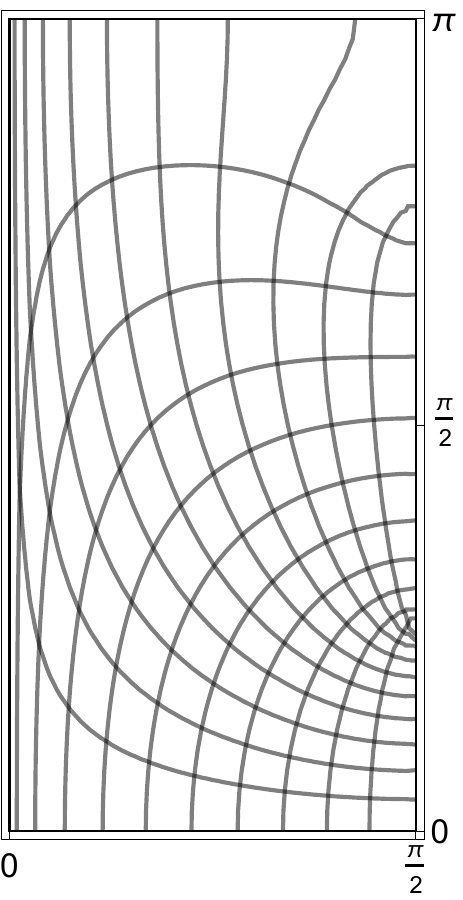}
  }
  \caption{Quarter Sphere Construction on the parameter plane. (a) Mesh with geometric grading to all singularities, including
  a strong geometric refinement along the edge, corresponding to the singularity at the pole.
    (b), (c) Contour lines of the solutions of the original and conjugate problems. (d) Projection of the conformal map onto
  the parameter domain.}
  \label{fig:qsfem}
\end{figure}

In the Quarter Sphere problem one of the corners of the quadrilateral lies on the pole. This singularity becomes an edge singularity
due to mapping. The resulting map is shown in Figure~\ref{fig:qsgrid}. Therefore, it is not sufficient to refine only near the two singular points on the rim, but one has to add
strong refinement along the whole edge. This is illustrated in Figure~\ref{fig:qsfemA}. Other steps in the construction are as in the 
example above. 

\begin{figure}
  \centering
  \subfloat[Schwarz]{\label{fig:sphereconvergenceA}
    \includegraphics[width=0.45\textwidth]{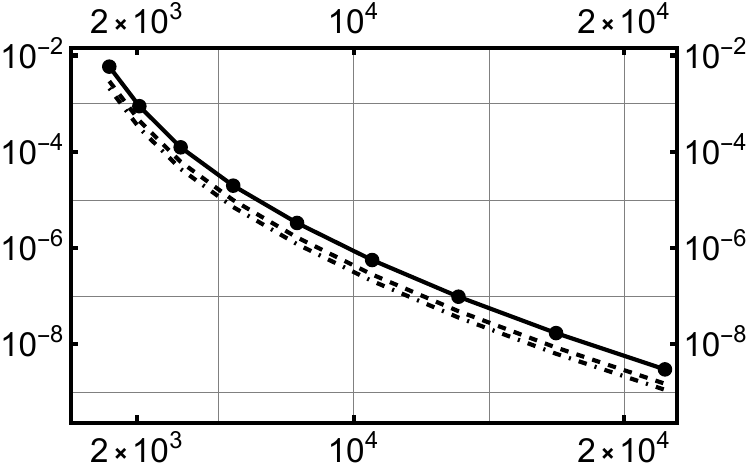}
  }
  \subfloat[Quarter Sphere]{\label{fig:sphereconvergenceB}
    \includegraphics[width=0.45\textwidth]{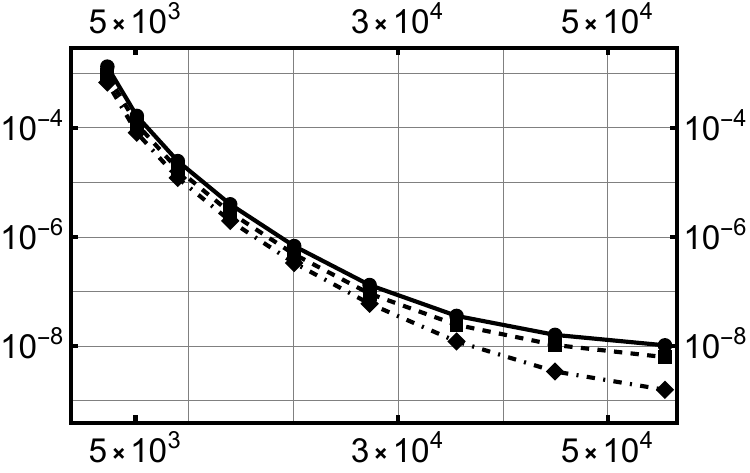}
  }
  \caption{Convergence plots. (a) Schwarzian hemisphere, (b) Quarter sphere. In both cases log-scale has been used.
  Three errors vs $N$ (the number of degrees of freedom) are shown: (i) reciprocal error (solid line), (ii) exact error in $M(Q)$ (dashed line), and
  (iii) estimated error in $M(Q)$ (dot-dashed line). Convergence is measured using a fixed mesh with a constant
  polynomial order $p=2,\ldots,10$. In (a) the observed, convergence is exponential over the whole range, in (b) there is a slight loss of rate for $p > 8$.
 }
  \label{fig:sphereconvergence}
\end{figure}

The four points $z_i$ on the boundary are
\[
z_1 = (0,0,1),\;
z_2=(1,0,0),\;
z_3=\Big(\frac{1}{\sqrt{2}},\;
\frac{1}{\sqrt{2}},0\Big),\;
z_4=\Big(-\frac{1}{\sqrt{2}},\frac{1}{\sqrt{2}},0\Big).
\]
In this case the exact moduli are not known. However, there is considerable computational evidence that
the following conjecture is true.
\begin{conjecture}
  $M(Q) = \sqrt{2}, \; M(\widetilde{Q}) = 1/\sqrt{2}$.
\end{conjecture}
The observed exponential convergence is shown in Figure~\ref{fig:sphereconvergenceB}.
The performance is not as good as in the case of the Schwarzian hemisphere, probably due to the edge singularity.
Notice, that the reciprocal error and the exact error are perfectly aligned, whereas the
auxiliary space error estimator is optimistic as the polynomial order increases.

\subsection{Cusps: Hyperbolic Quadrilateral on Sphere}
Let $Q_s$ be the quadrilateral whose sides are circular arcs
perpendicular to the unit circle with vertices $e^{is}$, $e^{(\pi-s)i}$,
$e^{(s-\pi)i}$ and $e^{-si}$. We call quadrilaterals of this type 
hyperbolic quadrilaterals as their sides are geodesics in the hyperbolic 
geometry of the unit disk \cite{hrv3}. For the general case we use notation
$\text{HypQuad}((x_0,y_0),R,s)$, where the unit circle is shifted to $(x_0,y_0)$ and
scaled to radius $R$ (For this particular instance, see Table~\ref{tbl:experiments}).

These domains are interesting since they have four cusps at the corners. The reciprocal relation
makes it possible to monitor the accuracy of the solutions even in such complicated domains.
If the configuration considered here would be constrained to the plane, then the computed
moduli would both be identical to one due to symmetry. Mapping onto the sphere
naturally perturbs this balance. We obtained the moduli
\begin{equation}
  M(Q) =1.8062303587451534, \; M(\widetilde{Q}) = 0.5536392383024755,
\end{equation}
with $\operatorname{reci}(Q) = 1.447375552743324\times 10^{-11}$.
\begin{figure}
  \centering
  \subfloat[]{\label{fig:hypq}
    \includegraphics[height=2.25in]{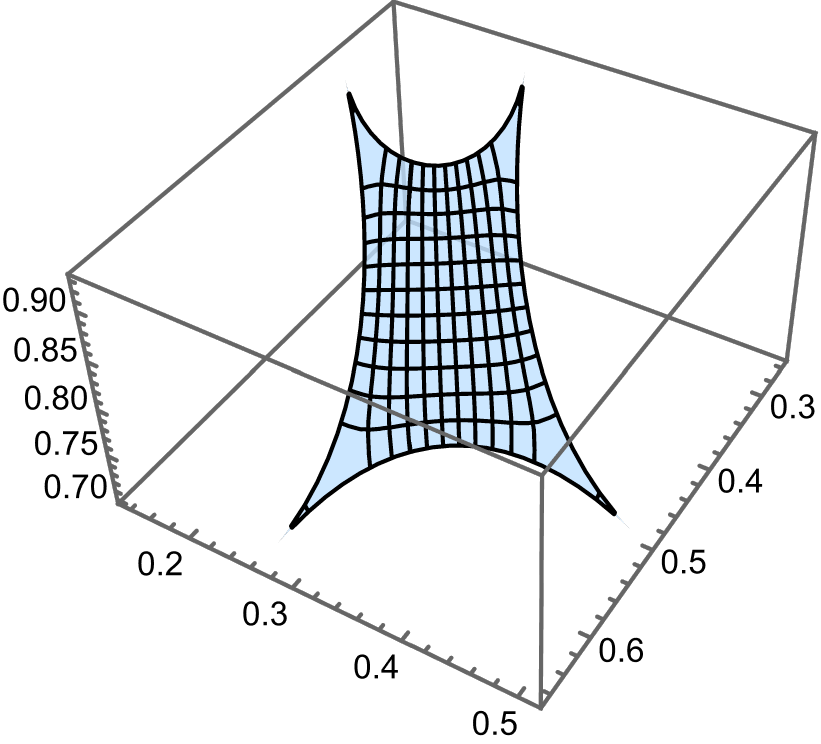}
  }
  \subfloat[]{\label{fig:mcs}
    \includegraphics[height=2.25in]{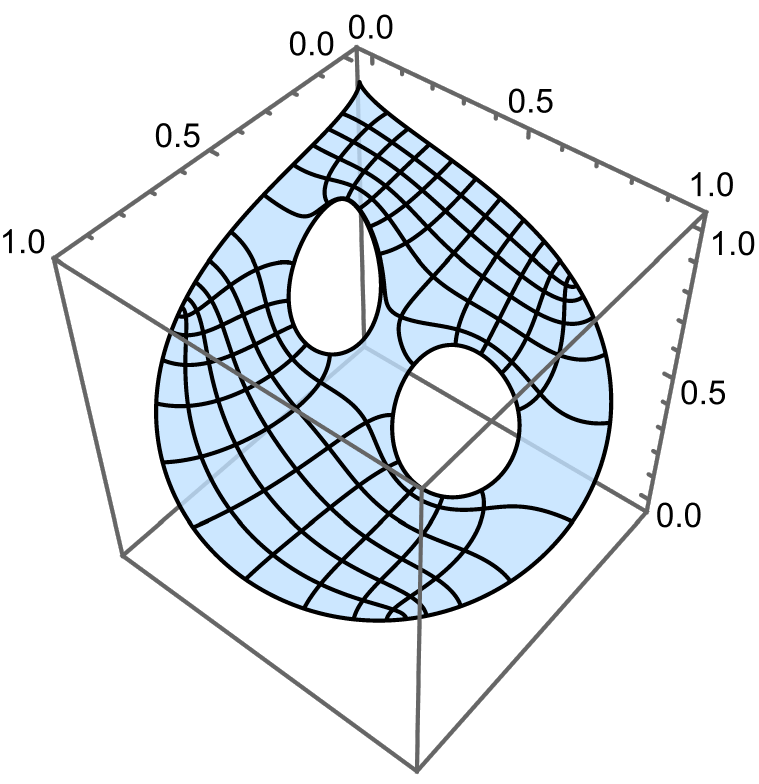}
  }
  \caption{Featured examples of maps. (a) Hyperbolic quadrilateral. (b) Multiply connected domain.}
  \label{fig:cusps}
\end{figure}

%

\subsection{Multiply Connected Domain: Two Holes}
\begin{figure}
  \centering
  \subfloat[]{\label{fig:holesA}
    \includegraphics[height=2in]{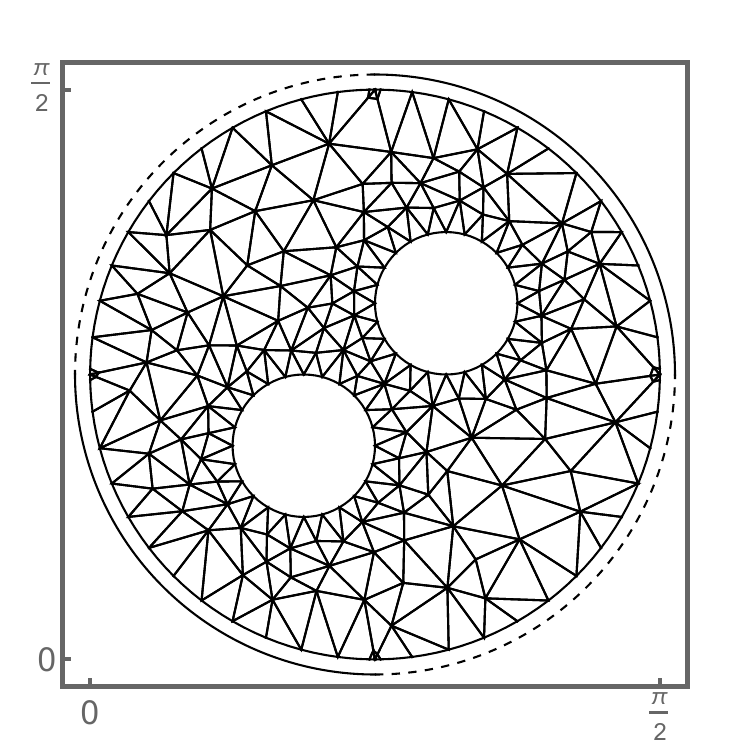}
  }
  \subfloat[]{\label{fig:holesB}
    \includegraphics[width=0.55\textwidth]{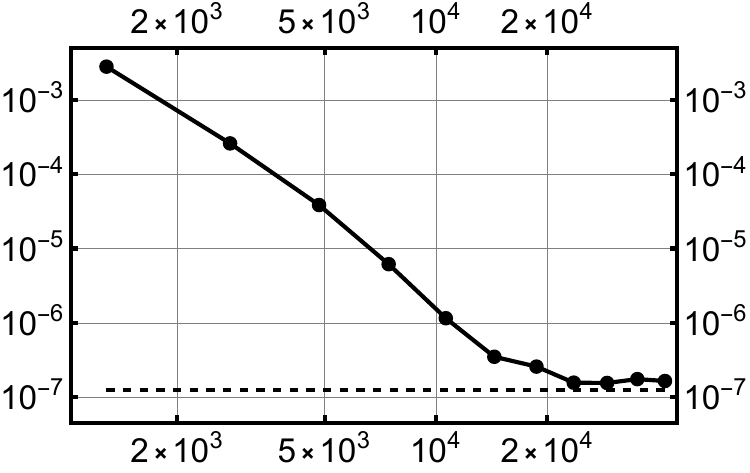}
  }
  \caption{Multiply Connected Domain. (a) Mesh in the parameter plane. (b) Convergence of the reciprocal error as a function of degrees of freedom. The threshold limit of $10^{-7}$ indicates the best possible accuracy given that the optimisation minimises the square of the error (loglog-plot).}
  \label{fig:holes}
\end{figure}

The conjugate function method can also be defined for multiply connected domains.
Of course, the immediate question concerns the definition and subsequent construction of the
conjugate problem. Here the so-called $Q$-type variant has been considered. The algorithm can be
defined as follows \cite{hqr2}:
\begin{description}
\item[Assumptions] The outer boundary has four corners which define a quadrilateral.
\item[Initialization] Solve the problem on the quadrilateral with zero Neumann boundary conditions on the holes.
\item[Optimization] Construct the conjugate problem on the outer boundary and find the potentials on the holes via optimization. The objective is to minimize the reciprocal error.
\end{description}
The mesh used is shown in Figure~\ref{fig:holesA}. Once the problem on the quadrilateral with
Neumann holes has been solved, the potentials on the holes in the conjugate problem are found
using standard optimization. Here the interior point method is applied as implemented in Mathematica.
The convergence rate is exponential up to the accuracy limit imposed by the chosen optimization
algorithm, see Figure~\ref{fig:holesB}.

The precise description of the problem on the parameter plane is as follows: 
As indicated in Table~\ref{tbl:experiments}, the computational domain is a unit circle centred at 
$(1/2,1/2)$, with two circular holes with radius $= 1/4$ and centres at $(1/4,1/4)$ and $(3/4,3/4)$
denoted as $B_1$ and $B_2$, respectively. The four corners are selected symmetrically at 
$z_i = (1/2 + \cos(\pi (k-1) /2), 1/2 + \sin(\pi (k-1) /2))$, $k=1,\dots,4$. 
The optimization step gives us potentials $v_1 = 0.5343446377370098$ on $\partial B_1$ and 
$v_2 = 0.5343446377370098$ on $\partial B_2$. Hence the two problems can be summarized as:
\begin{equation}
  \begin{cases}
    \Delta_S u = 0,\quad \text{ in } \Omega, \\
    u = 0, \quad \text{ on } \gamma_1,  \\
    u = 1, \quad \text{ on } \gamma_3, \\
    \partial u/\partial n = 0, \quad  \text{ on } \gamma_2, \gamma_4, \\
    \partial u/\partial n = 0, \quad  \text{ on } \partial B_1, \partial B_2, \\
  \end{cases} \text{ leading to }
    \begin{cases}
    \Delta_S v = 0,\quad \text{ in } \Omega, \\
    v = 0, \quad \text{ on } \gamma_2, \\ 
    v = 1, \quad \text{ on } \gamma_4, \\ 
    v = v_1, \quad \text{ on } \partial B_1,\\
    v = v_2, \quad \text{ on } \partial B_2,\\
    \partial v/\partial n = 0, \quad  \text{ on } \gamma_1, \gamma_3. 
  \end{cases}
\end{equation}
We obtained moduli
\begin{equation}
  M(Q) = 0.7901907571620941, \; M(\widetilde{Q}) = 1.2655174148067712,
\end{equation}
with $\operatorname{reci}(Q) = 1.6420797832594758\times 10^{-7}$.

\subsection{Seashell}
\begin{figure}
  \centering
  \subfloat{
    \includegraphics[height=3in]{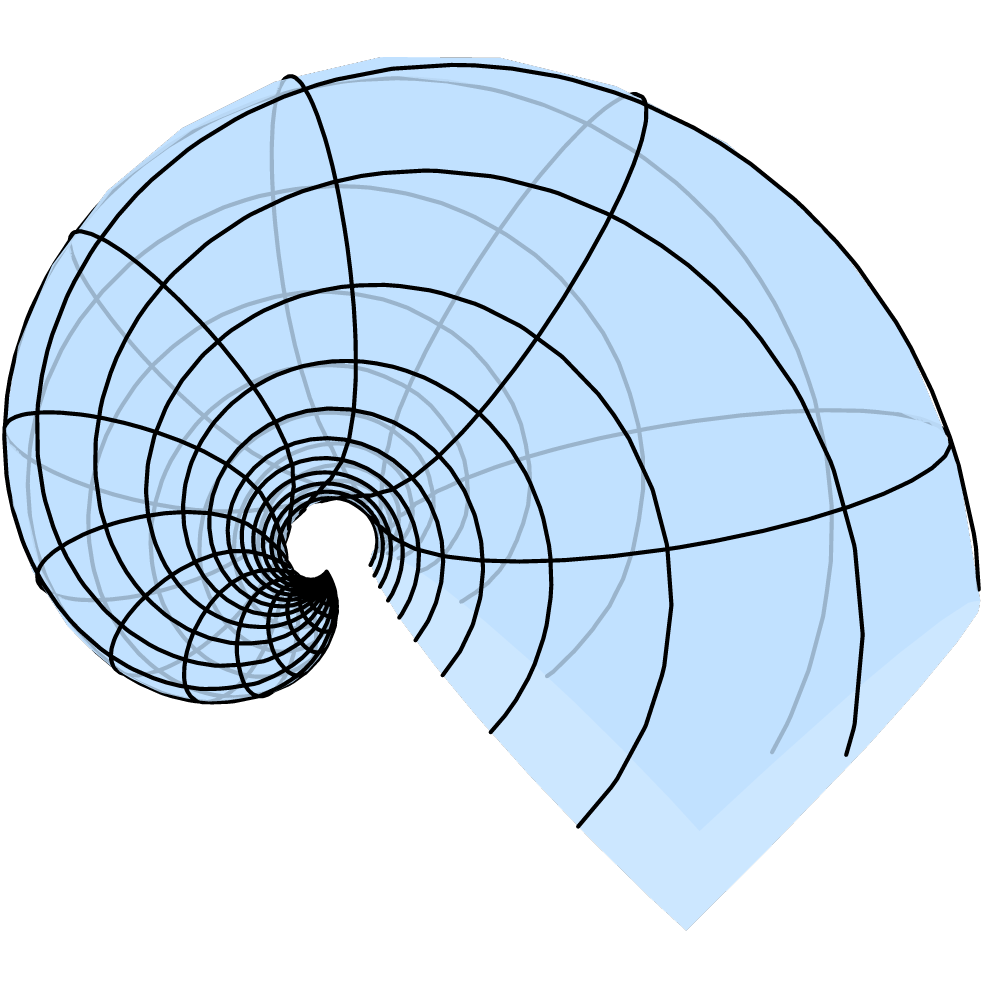}
  }
  \caption{Seashell. Example of a non-smooth surface.}
  \label{fig:seashell}
\end{figure}

Our last example is a so-called seashell surface
\begin{equation}\label{seashell-parametric}
\left\{\begin{array}{rcl}
x &=&a(1-v/(2\pi))\cos(n v)(1+\cos(u))+c\cos(n v), \\
y &=&a(1-v/(2\pi))\sin(n v)(1+\cos(u))+c\sin(n v), \\ 
z &=&b v/(2\pi)+a(1-v/(2\pi))\sin(u),
\end{array}\right.
\end{equation}
where we let $n=1$, $a=1$, $b=1$, and $c=1/10$, with $u \in [0,2\pi]$, $v \in [-\pi,\pi]$.
The resulting map is illustrated in Figure~\ref{fig:seashell}.

\begin{figure}
  \centering
  \subfloat[]{\label{fig:seashellDetailA}
    \includegraphics[width=0.45\textwidth]{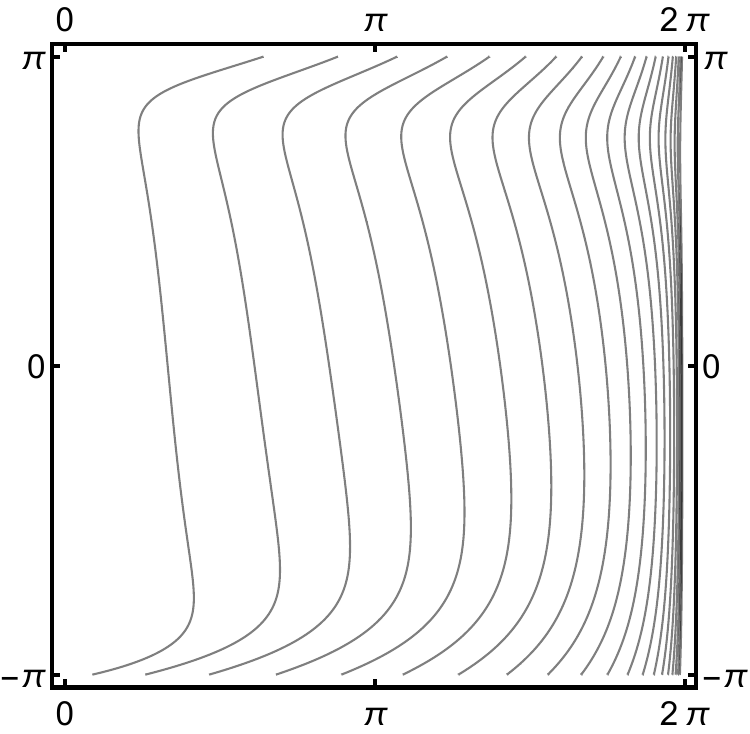}
  }
  \subfloat[]{\label{fig:seashellDetailB}
    \includegraphics[width=0.45\textwidth]{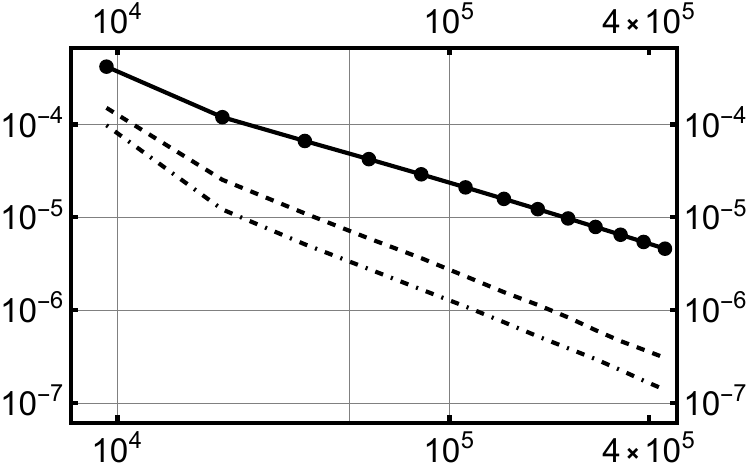}
  }
  \caption{Seashell. (a) Contour lines on $\widetilde{Q}$ where the effect of the "tip" is clearly visible. Mapping onto the surface connects 
  the contour lines periodically. (b) Three errors vs $N$ (the number of degrees of freedom) are shown: (i) reciprocal error (solid line), (ii) estimated error in $M(Q)$ (dashed line), and
  (iii) estimated error in $M(\widetilde{Q})$ (dot-dashed line), (loglog-plot). }
  \label{fig:seashellDetail}
\end{figure}

This case is challenging, since the surface is not smooth -- for details on this class of problems we refer to Miklyukov \cite{Miklyukov}. 
The parameterization is singular along the edge $u = 2\pi$, and hence we choose an offset $\epsilon = 1/10$
as in the case of Section~\ref{sec:earth}. There are also strong boundary layer -like features along the edges $v = \pm \pi$, and
therefore the simplest mesh is a regular grid with strong grading along the three boundaries requiring it.
The effect of the sharp tip is visible in the contour lines (see Figure~\ref{fig:seashellDetailA})suggesting that strong mesh refinement is necessary at that part of the domain.
With proper grading the convergence rate is slightly better than algebraic (see Figure~\ref{fig:seashellDetailB}).
We obtained the moduli
\begin{equation}
  M(Q) = 1.567020274702868, \; M(\widetilde{Q}) = 0.638156772214456,
\end{equation}
with $\operatorname{reci}(Q) = 4.600498992424207\times 10^{-6}$.

\subsection{On Computational Complexity}
\label{sec:computational_complexity}
The proposed method reduces to a 2D $hp$-FEM with variable coefficients over the parameter domain.
The computational complexity is, as usual with $hp$-FEM, dominated by the numerical integration rather
than the solution of the linear systems. For details on the quadrature techniques used, 
we refer to \cite{Hakula2013}. For linear systems, Cholesky factorization is used in all cases.

Let us revisit the seashell example. In the final mesh there are 2365 nodes, 4632 edges, 
and 2268 quadrilateral elements. With a constant $p = 14$ the number of degrees of freedom is 445873, and
the dimension of the auxiliary error space is 131640. The elemental quadrature rules are standard
Gauss rules with $(p+q)^2$ points, where the extension $q = 5$. 
On average, after the integration of the full system in 2700 seconds, 
the subsequent solves for $p=2,\ldots,14$ with the error estimation taking 280 seconds, 
using one Mathematica 14.0 kernel
on an Apple Silicon Macbook Pro (M3 Max, 2023 model).

\section{Conclusions}
\label{sec:conclusions}
We have the conjugate function method for computation of conformal
mappings on surfaces. Both theoretical and implementation aspects have been addressed
and the efficacy of the method has been demonstrated with a series of numerical experiments
covering problems with specific features such as singularities and cusps, including
multiply connected domains. The
fundamental ideas can be applied to address a much wider class of problems, however.
Our method relies on numerical solution of PDEs and can be applied to
general surfaces using standard tools as long as the suitable global parameterization is available. 
In terms of computational cost, the method is
competitive especially in cases where the $p$-version of FEM can be applied directly. 

\subsection*{Acknowledgments}
The second author (A. Rasila) was partly supported by NSF of China under the number 11971124, NSF of Guangdong Province under the numbers 2021A1515010326 and 2024A1515010467, and Li Ka Shing Foundation.

\bibliographystyle{siamplain}
\bibliography{references}
\end{document}